\documentclass[11pt]{amsart}
\usepackage{amscd,amssymb}
\usepackage{amsthm,amsmath,amssymb}
\usepackage[matrix,arrow]{xy}

\newtheorem{theorem}[equation]{Theorem}

\newtheorem{lemma}[equation]{Lemma}
\newtheorem{corollary}[equation]{Corollary}

\theoremstyle{definition}
\newtheorem{example}[equation]{Example}

\theoremstyle{remark}
\newtheorem{remark}[equation]{Remark}

\makeatletter\@addtoreset{equation}{section} \makeatother

\author{Ivan Cheltsov}

\title{Factorial threefold hypersurfaces}

\address{\begin{tabbing}
\hspace*{28 em}\=\kill
School of Mathematics\\
University of Edinburgh\\
Edinburgh EH9 3JZ, UK\\
\\
\texttt{I.Cheltsov@ed.ac.uk}
\end{tabbing}}

\thanks{We assume that all varieties are projective, normal, and defined over $\mathbb{C}$.}%

\begin{document}

\begin{abstract}
Let $X$ be a hypersurface in $\mathbb{P}^{4}$ of degree~$d$ that
has at most isolated ordinary double points. We prove that $X$ is
factorial in the case when $X$ has at most $(d-1)^{2}-1$ singular
points.
\end{abstract}

\maketitle

\section{Introduction}
\label{section:intro}

The Cayley--Bacharach theorem (see \cite{DaGeOr85}, \cite{EiJ96}),
in its classical form, may be seen as~a~result about the number of
independent linear conditions imposed on forms of a given degree
by a~certain finite subset of $\mathbb{P}^{n}$. The purpose of
this paper is to prove the following result.

\begin{theorem}
\label{theorem:main} Let $\Sigma$ be a finite subset in
$\mathbb{P}^{n}$, and let $\mu$ be a natural number such that
\begin{itemize}
\item  the inequalities $\mu\geqslant 2$ and $|\Sigma|\leqslant\mu^{2}-1$ hold,%
\item at most $\mu k$ points in the set $\Sigma$ lie on a~curve in $\mathbb{P}^{n}$ of degree~$k=1,\ldots,\mu-1$,%
\end{itemize}
where $n\geqslant 2$. Then $\Sigma$ imposes independent linear
conditions on forms~of~degree~$2\mu-3$.
\end{theorem}

Let $X$ be a hypersurface in $\mathbb{P}^{4}$ of
degree~$d\geqslant 3$ such that the threefold $X$ has at most
isolated ordinary double points. Then $X$ can be given by the
equation
$$
f\big(x,y,z,t,u\big)=0\subset\mathbb{P}^{4}\cong\mathrm{Proj}\Big(\mathbb{C}\big[x,y,z,t,u\big]\Big),%
$$
where $f(x,y,z,t,u)$ is a homogeneous polynomial of degree $d$.

\begin{remark}
\label{remark:factoriality} It follows from \cite{Ha70}
and~\cite{Di90} that the following conditions are equivalent:
\begin{itemize}
\item every Weil divisor on the threefold $X$ is a Cartier divisor;%
\item every surface $S\subset X$ is cut out on $X$ by a hypersurface in $\mathbb{P}^{4}$;%
\item the ring
$$
\mathbb{C}\big[x,y,z,t,u\big]\Big\slash\big\langle f\big(x,y,z,t,u\big)\big\rangle%
$$
is a unique factorization domain;%
\item the set $\mathrm{Sing}(X)$ imposes independent linear conditions on forms of degree $2d-5$.%
\end{itemize}
\end{remark}

We say that $X$ is factorial if every Weil divisor on $X$ is a
Cartier divisor.

\begin{example}
\label{example:non-factorial-hypersurface} Suppose that $X$ is
given by
$$
xg\big(x,y,z,t,u\big)+yh\big(x,y,z,t,u\big)=0\subset\mathbb{P}^{4}\cong\mathrm{Proj}\Big(\mathbb{C}\big[x,y,z,t,u\big]\Big),%
$$
where $g$ and $h$ are general homogeneous polynomials of degree
$d-1$.~Then
\begin{itemize}
\item the threefold $X$ has at most isolated ordinary double points,%
\item the equality $|\mathrm{Sing}(X)|=(d-1)^{2}$ holds, but $X$ is not fac\-to\-rial.%
\end{itemize}
\end{example}

The assertion of Theorem~\ref{theorem:main} implies the following
result (cf. \cite{CiGe05d}, \cite{Ch04t}, \cite{Ch08}).

\begin{theorem}
\label{theorem:Cheltsov} Suppose that
$|\mathrm{Sing}(X)|<(d-1)^{2}$. Then $X$ is factorial.
\end{theorem}

\begin{proof}
The set $\mathrm{Sing}(X)$ is a set-theoretic intersection of
hypersurfaces of degree $d-1$. Then
\begin{itemize}
\item  the inequalities $d-1\geqslant 2$ and $|\,\mathrm{Sing}(X)|\leqslant(d-1)^{2}-1$ hold,%
\item at most $(n-1)k$ points in the set $\mathrm{Sing}(X)$ lie on a~curve in $\mathbb{P}^{4}$ of degree~$k=1,\ldots,n-2$,%
\end{itemize}
which immediately implies that the points of the set
$\mathrm{Sing}(X)$ imposes independent linear conditions on forms
of degree $2d-5$ by Theorem~\ref{theorem:main}. Thus, the
threefold $X$ is factorial.
\end{proof}

The assertion of Theorem~\ref{theorem:Cheltsov} is proved in
\cite{Ch04e} and \cite{ChPa05d} in the case when $d\leqslant 7$.

\begin{remark}
\label{remark:Mella} Suppose that $d=4$ and $X$ is factorial. Then
it follows from \cite{Me03} that
\begin{itemize}
\item the threefold $X$ is non-rational,%
\item the threefold $X$ is not birational to a conic bundle,%
\item the threefold $X$ is not birational to a fibration into rational surfaces,%
\end{itemize}
but ge\-ne\-ral determinantal quartic hypersurfaces in
$\mathbb{P}^{4}$ are rational.
\end{remark}

The author thanks J.\,Park, Yu.\,Pro\-kho\-rov, V.\,Sho\-ku\-rov,
K.\,Shramov for useful comments.

\section{The proof}
\label{section:main}

Let $\Sigma$ be a finite subset in $\mathbb{P}^{n}$, and let $\mu$
be a natural number such that
\begin{itemize}
\item  the inequalities $\mu\geqslant 2$ and $|\Sigma|\leqslant\mu^{2}-1$ hold,%
\item at most $\mu k$ points in the set $\Sigma$ lie on a~curve in $\mathbb{P}^{n}$ of degree~$k=1,\ldots,\mu-1$,%
\end{itemize}
where $n\geqslant 2$. Suppose that $\Sigma$ imposes dependent
linear conditions on forms of degree $2\mu-3$.

\begin{remark}
\label{remark:mu-2} The inequality $\mu\geqslant 3$ holds.
\end{remark}

The following result is proved in \cite{Bes83} and \cite{DaGe88}.

\begin{theorem}
\label{theorem:Bese} Let $P_{1},\ldots,
P_{\delta}\in\mathbb{P}^{2}$ be distinct points such that
\begin{itemize}
\item at most $k(\xi+3-k)-2$ points in $\{P_{1},\ldots, P_{\delta}\}$ lie on a curve of degree $k\leqslant (\xi+3)/2$,%
\item the inequality
$$
\delta\leqslant
\mathrm{max}\left\{\Big\lfloor\frac{\xi+3}{2}\Big\rfloor\left(\xi+3-\Big\lfloor\frac{\xi+3}{2}\Big\rfloor\right)-1,
\Big\lfloor\frac{\xi+3}{2}\Big\rfloor^{2}\right\}%
$$
holds, where $\xi$ is a natural number such that $\xi\geqslant 3$,
\end{itemize}
and let $\pi\colon Y\to\mathbb{P}^2$ be a blow up of the points
$P_{1},\ldots, P_{\delta}$. Then the linear system
$$
\Big|\pi^{*}\Big(\mathcal{O}_{\mathbb{P}^2}\big(\xi\big)\Big)-\sum_{i=1}^{\delta}E_{i}\Big|
$$
does not have base points, where $E_{i}$ is the
$\pi$-ex\-cep\-ti\-onal divisor such that $\pi(E_{i})=P_{i}$.
\end{theorem}

There is a point $P\in\Sigma$ such that every
hypersurface\footnote{For simplicity we consider homogeneous forms
on $\mathbb{P}^{n}$ as hypersurfaces.} in $\mathbb{P}^{n}$ of
degree $2\mu-3$ that contains the set $\Sigma\setminus P$ must
contain~the point $P\in\Sigma$. Let us derive a contradiction.

\begin{lemma}
\label{lemma:surfaces-I} The inequality $n\ne 2$ holds.
\end{lemma}

\begin{proof}
Suppose that $n=2$. Let us prove that at most $k(2\mu-k)-2$ points
in $\Sigma\setminus P$ can lie on a curve of degree $k\leqslant
\mu$. It is enough to show that
$$
k\big(2\mu-k\big)-2\geqslant k\mu
$$
for every $k\leqslant \mu$. We must prove this only for
$k\geqslant 1$ such that
$$
k\big(2\mu-k\big)-2<\big|\Sigma\setminus P\big|\leqslant \mu^{2}-2,%
$$
because otherwise the condition that at most $k(2\mu-k)-2$ points
in the set $\Sigma\setminus P$ can lie on a curve of degree $k$ is
vacuous. Therefore, we may assume that $k<\mu$.

We may assume that $k\ne 1$, because at most $\mu\leqslant 2\mu-3$
points of $\Sigma\setminus P$ lie on a line. Then
$$
k\big(2\mu-k\big)-2\geqslant k\mu\iff \mu>k,
$$
which implies that at most $k(2\mu-k)-2$ points in
$\Sigma\setminus P$ can lie on a curve in $\mathbb{P}^{2}$ of
degree $k\leqslant \mu$.

Thus, it  follows from Theorem~\ref{theorem:Bese} that there is a
curve of degree $2\mu-3$ that contains all points of the set
$\Sigma\setminus P$ and does not contain the point $P\in\Sigma$,
which is a contradiction.
\end{proof}

Moreover, we may assume that $n=3$ due to the following result.

\begin{lemma}
\label{lemma:zero-dimensional} Let $\Lambda\subset\Sigma$ be a
subset, let $\psi\colon\mathbb{P}^n\dasharrow\mathbb{P}^m$ be a
general linear projection, and let
$$
\mathcal{M}\subseteq\Big|\mathcal{O}_{\mathbb{P}^{n}}\big(k\big)\Big|
$$
be a linear subsystem~that contains all hypersurfaces that pass
through $\Lambda$. Suppose that
\begin{itemize}
\item the inequality and $|\Lambda|\geqslant\mu k+1$ holds,%
\item the set $\psi(\Lambda)$ is contained in an irre\-du\-ci\-ble reduced curve of degree $k$,%
\end{itemize}
where $n>m\geqslant 2$. Then $\mathcal{M}$ has no base curves, and either $m=2$,~or~$k>\mu$.%
\end{lemma}

\begin{proof}
We may assume that there are linear subspaces $\Omega$ and
$\Pi\subset\mathbb{P}^n$ such that
$$
\psi\colon\mathbb{P}^n\dasharrow\Pi\cong\mathbb{P}^m
$$
is a projection from $\Omega$, where $\mathrm{dim}(\Omega)=n-m-1$
and $\mathrm{dim}(\Pi)=m$.

Suppose that there is an irreducible curve
$Z\subset\mathbb{P}^{n}$ such that $Z$ is contained in the base
locus of the~linear system $\mathcal{M}$. Put $\Xi=Z\cap\Lambda$.
We may assume that $\psi\vert_{Z}$ is a birational morphism,~and
$$
\psi\big(Z\big)\cap\psi\Big(\Lambda\setminus\Xi\Big)=\varnothing,
$$
because the projection $\psi$ is general. Then
$\mathrm{deg}(\psi(Z))=\mathrm{deg}(Z)$.

Let $C\subset\Pi$ be an irreducible curve of degree $k$ that
contains $\psi(\Lambda)$, and let $W\subset\mathbb{P}^{n}$ be
the~cone over the curve $C$ whose vertex is $\Omega$. Then $W\in
\mathcal{M}$, which implies that $Z\subset W$. We have
$$
\psi\big(Z\big)=C,
$$
which immediately implies that $\Xi=\Lambda$ and
$\mathrm{deg}(Z)=k$. But $|Z\cap\Sigma|\leqslant\mu k$, which is a
contradiction. Therefore, the linear system $\mathcal{M}$ does not
have base curves.

Now we suppose that $m\geqslant 3$ and $k\leqslant \mu$. Let us
show that this assumption leads to a~contradiction. Without loss
of generality, we may assume that $m=3$ and $n=4$.

Let $\mathcal{Y}$ be the set of all irreducible reduced surfaces
in $\mathbb{P}^{4}$ of degree $k$ that contains the set~$\Lambda$,
and let $\Upsilon$ be a subset of $\mathbb{P}^{4}$ that consists
of all points that are contained in every surface of the~set
$\mathcal{Y}$. Then $\Lambda\subseteq\Upsilon$. Arguing as above,
we see that $\Upsilon$ is a finite set.

Let $\mathcal{S}$ be the set of all surfaces in $\mathbb{P}^{3}$
of degree $k$ such that
$$
S\in\mathcal{S}\iff \exists\ Y\in\mathcal{Y}\ \text{such that}\ \psi\big(Y\big)=S\ \text{and}\ \psi\big\vert_{Y}\ \text{is a birational morphism},%
$$
and let $\Psi\subset\mathbb{P}^{3}$ that  consists of all points
contained in every surface in $\mathcal{S}$. Then
$\mathcal{S}\ne\varnothing$ and
$$
\psi\big(\Lambda\big)\subseteq\psi\big(\Upsilon\big)\subseteq\Psi.
$$

For every point $O\in\Pi\setminus\Psi$ and for a general surface
$Y\in\mathcal{Y}$, we may assume that the line passing through $O$
and $\Omega$ does not intersect $Y$. But $\psi\vert_{Y}$ is a
birational morphism. Then
$$
\psi\big(\Upsilon\big)=\Psi,
$$
and  $\psi(\Lambda)\subseteq\Psi$ contains at least $\mu
k+1\geqslant k^{2}+1$ points that are contained in a~curve of
degree~$k$, which is impossible, because $\Psi$ is a set-theoretic
intersection of surfaces of degree $k$.
\end{proof}

Fix a sufficiently general hyperplane $\Pi\subset\mathbb{P}^{3}$.
Let
$$
\psi\colon \mathbb{P}^{3}\dasharrow\Pi\cong\mathbb{P}^{2}
$$
be a projection from a
sufficiently general point $O\in\mathbb{P}^{3}$. Put
$\Sigma^{\prime}=\psi(\Sigma)$ and $P^{\prime}=\psi(P)$.

\begin{lemma}
\label{lemma:surfaces-II} There is a curve $C\subset\Pi$ of degree
$k\leqslant\mu-1$ such that $|C\cap\Sigma^{\prime}|\geqslant \mu
k+1$.
\end{lemma}

\begin{proof}
We  suppose that at most $\mu k$ points of the set
$\Sigma^{\prime}$ are contained in a curve in $\Pi$ of degree~$k$
for every $k\leqslant \mu-1$. Then arguing as in the proof of
Lemma~\ref{lemma:surfaces-I}, we obtain a curve
$$
Z\subset\Pi\cong\mathbb{P}^{2}
$$
of degree $2\mu-3$ that contains the set $\Sigma^{\prime}\setminus
P^{\prime}$ and does not pass through the point $P^{\prime}$.

Let $Y$ be the cone in $\mathbb{P}^{3}$ over the curve $Z$ whose
vertex is the point $O$. Then $Y$ is a~surface of degree $2\mu-3$
that contains all points of the set $\Sigma\setminus P$ but does
not contain the point $P\in\Sigma$.
\end{proof}

It immediately follows from Lemma~\ref{lemma:zero-dimensional}
that $k\geqslant 2$.

\begin{lemma}
\label{lemma:conic} Suppose that $|C\cap\Sigma^{\prime}|\geqslant
9$. Then $k\geqslant 3$.
\end{lemma}

\begin{proof}
Suppose that $k=2$. Let $\Phi\subseteq\Sigma$ be a subset such
that $|\Phi|\geqslant 9$, but $\psi(\Phi)$ is contained in the
conic $C\subset\Pi$. Then the conic $C$ is irreducible by
Lemma~\ref{lemma:zero-dimensional}.

Let $\mathcal{D}$ be a linear system of quadric hypersurfaces in
$\mathbb{P}^{3}$ containing $\Phi$. Then $\mathcal{D}$ does not
have base curves by Lemma~\ref{lemma:zero-dimensional}. Let $W$ be
a cone in $\mathbb{P}^{3}$ over $C$ with the vertex $\Omega$. Then
$$
8=D_{1}\cdot D_{2}\cdot W\geqslant
\sum_{\omega\in\Phi}\mathrm{mult}_{\omega}(D_{1})\mathrm{mult}_{\omega}(D_{2})\geqslant |\,\Phi|\geqslant 9,%
$$
where $D_{1}$ and $D_{2}$ are general divisors in the linear
system $\mathcal{D}$.
\end{proof}

We may assume that $k$ is the smallest natural number such that at
least $\mu k+1$ points in $\Sigma^{\prime}$ lie on a curve of
degree $k$. Then there is a non-empty disjoint union
$$
\bigcup_{j=k}^{l}\bigcup_{i=1}^{c_{j}}\Lambda_{j}^{i}\subset\Sigma
$$
such that $|\Lambda_{j}^{i}|\geqslant \mu j+1$, all points of the
the set $\psi(\Lambda_{j}^{i})$ are contained in an irreducible
reduced curve of degree~$j$, and at most $\mu\zeta$ points of the
subset
$$
\psi\left(\Sigma\setminus\Big(\bigcup_{j=k}^{l}\bigcup_{i=1}^{c_{j}}\Lambda_{j}^{i}\Big)\right)\subsetneq\Sigma^{\prime}\subset\Pi\cong\mathbb{P}^{2}
$$
lie on a curve in $\Pi$ of degree $\zeta$ for every natural number
$\zeta$. Put
$$
\Lambda=\bigcup_{j=k}^{l}\bigcup_{i=1}^{c_{j}}\Lambda_{j}^{i}.
$$

Let $\Xi_{j}^{i}$ be the base locus of the linear subsystem in
$|\mathcal{O}_{\mathbb{P}^{3}}(j)|$ that contains all surfaces
passing through the set $\Lambda_{j}^{i}$. Then $\Xi_{j}^{i}$ is a
finite set by Lemma~\ref{lemma:zero-dimensional}, and
\begin{equation}
\label{equation:number-of-good-points}
\big|\Sigma\setminus\Lambda\big|\leqslant\mu\left(\mu-\sum_{i=k}^{l}c_{i}\mu i\right)-2.%
\end{equation}

\begin{corollary}
\label{corollary:from-the-number-of-good-points} The inequality
$\sum_{i=k}^{l}ic_{i}\leqslant \mu-1$ holds.
\end{corollary}

Put
$\Delta=\Sigma\cap(\cup_{j=k}^{l}\cup_{i=1}^{c_{j}}\Xi_{j}^{i})$.
Then $\Lambda\subseteq\Delta\subseteq\Sigma$.

\begin{lemma}
\label{lemma:non-vanishing} The set $\Delta$ impose independent
linear conditions on forms of degree $2\mu-3$.
\end{lemma}

\begin{proof}
Let us consider the subset $\Delta\subset\mathbb{P}^{3}$ as a
closed subscheme of $\mathbb{P}^{3}$, and let
$\mathcal{I}_{\Delta}$ be the ideal sheaf of the subscheme
$\Delta$. Then there is an exact sequence
$$
0\longrightarrow\mathcal{I}_{\Delta}\otimes\mathcal{O}_{\mathbb{P}^{3}}\big(2\mu-3\big)\longrightarrow\mathcal{O}_{\mathbb{P}^{3}}\big(2\mu-3\big)\longrightarrow\mathcal{O}_{\Delta}\longrightarrow 0,%
$$
which implies that $\Delta$ imposes independent conditions on
forms of degree $2\mu-3$ if and only if
$$
h^{1}\Big(\mathcal{I}_{\Delta}\otimes\mathcal{O}_{\mathbb{P}^{3}}\big(2\mu-3\big)\Big)=0.
$$

Suppose
$h^{1}(\mathcal{I}_{\Delta}\otimes\mathcal{O}_{\mathbb{P}^{3}}(2\mu-3))\ne
1$. Let us show that this assumption leads to a contradiction.

Let $\mathcal{M}$ be a linear subsystem in
$|\mathcal{O}_{\mathbb{P}^{3}}(\mu-1)|$ that contains all surfaces
that pass through all point of the set $\Delta$. Then the base
locus of $\mathcal{M}$ is zero-dimensional, because
$\sum_{i=k}^{l}ic_{i}\leqslant \mu-1$ and
$$
\Delta\subseteq\bigcup_{j=k}^{l}\bigcup_{i=1}^{c_{j}}\Xi_{j}^{i},
$$
but $\Xi_{j}^{i}$ is a zero-dimensional base locus of a linear
subsystem in $|\mathcal{O}_{\mathbb{P}^{3}}(j)|$. Put
$$
\Gamma=M_{1}\cdot M_{2}\cdot M_{3},
$$
where $M_{1},M_{2},M_{3}$ are general surfaces in the linear
system $\mathcal{M}$. Then $\Gamma$ is a zero-dimensional
subscheme of $\mathbb{P}^{3}$, and $\Delta$ is a closed subscheme
of the scheme $\Gamma$.

Let $\Upsilon$ be a closed subscheme of the scheme $\Gamma$ such
that
$$
\mathcal{I}_{\Upsilon}=\mathrm{Ann}\Big(\mathcal{I}_{\Delta}\big\slash\mathcal{I}_{\Gamma}\Big),
$$
where $\mathcal{I}_{\Upsilon}$  and $\mathcal{I}_{\Gamma}$ are the
ideal sheaves of the subschemes $\Upsilon$  and $\Gamma$,
respectively. Then
$$
0\ne
h^{1}\Big(\mathcal{O}_{\mathbb{P}^{3}}\big(2\mu-3\big)\otimes\mathcal{I}_{\Delta}\Big)=h^{0}\Big(\mathcal{O}_{\mathbb{P}^{3}}\big(\mu-4\big)\otimes\mathcal{I}_{\Upsilon}\Big)-h^{0}\Big(\mathcal{O}_{\mathbb{P}^{3}}\big(\mu-4\big)\otimes\mathcal{I}_{\Gamma}\Big)%
$$
by Theorem~3 in \cite{DaGeOr85} (see also \cite{EiJ96}). Thus,
there is a surface
$F\in|\mathcal{O}_{\mathbb{P}^{n}}(\mu-4)\otimes\mathcal{I}_{\Upsilon}|$.
Then
$$
\big(\mu-4\big)\big(\mu-1\big)^{2}=F\cdot M_{1}\cdot M_{2}\geqslant h^{0}\big(\mathcal{O}_{\Upsilon}\big)=h^{0}\big(\mathcal{O}_{\Gamma}\big)-h^{0}\big(\mathcal{O}_{\Delta}\big)=\big(\mu-1\big)^{3}-\big|\Delta\big|,%
$$
which implies that $|\Delta|\geqslant 3(\mu-1)^{2}$. But
$|\Delta|\leqslant|\Sigma|<\mu^{2}$, which is impossible, because
$\mu\geqslant 3$.
\end{proof}

We see that $\Delta\subsetneq\Sigma$. Put
$\Gamma=\Sigma\setminus\Delta$ and
$d=2\mu-3-\sum_{i=k}^{l}ic_{i}$.

\begin{lemma}
\label{lemma:swapping} The set $\Delta$ imposes dependent linear
conditions on forms of degree $d$.
\end{lemma}

\begin{proof}
Suppose that the points of the set $\Delta$ impose independent
linear conditions on homogeneous polynomials of degree $d$. Let us
show that this assumption leads to a contradiction.

The construction of $\Delta$ implies the existence of a
homogeneous form $H$ of degree $\sum_{i=k}^{l}ic_{i}$~that
vanishes at all points of the set $\Delta$ and does not vanish at
any point of the set $\Gamma$.

Suppose that $P\in\Delta$. Then there is a homogenous form $F$ of
degree $2\mu-3$ that vanishes at every point of the set
$\Delta\setminus P$ and does not vanish at the point $P$ by
Lemma~\ref{lemma:non-vanishing}. Put
$$
\Gamma=\Big\{Q_{1},\ldots,Q_{\gamma}\Big\},
$$
where $Q_{i}$ is a point in $\Gamma$. Then
there~is~a~homo\-ge\-neous form $G_{i}$ of degree $d$ that
vanishes at every point in $\Gamma\setminus Q_{i}$ and does not
vanish at the point $Q_{i}$. Then
$$
F\big(Q_{i}\big)+\mu_{i}HG_{i}\big(Q_{i}\big)=0
$$
for some $\mu_{i}\in\mathbb{C}$, because $G_{i}(Q_{i})\ne 0$. Then
the homogenous form
$$
F+\sum_{i=1}^{\gamma}\mu_{i}HG_{i}
$$
vanishes on the set $\Sigma\setminus P$ and does not vanish at the
point $P$, which is a contradiction.

We see that $P\in\Gamma$. Then there is a~homo\-ge\-neous form $G$
of degree $d$ that vanishes at every point in $\Gamma\setminus P$
and does not vanish at $P$. Then $HG$ vanishes at every point of
the set $\Sigma\setminus P$ and does not vanish at the point $P$,
which is a contradiction.
\end{proof}

Put $\Gamma^{\prime}=\psi(\Gamma)$. Let us check that
$\Gamma^{\prime}$ and $d$ satisfy the hypotheses of
Theorem~\ref{theorem:Bese}.

\begin{lemma}
\label{lemma:d-3} The inequality $d\geqslant 3$ holds.
\end{lemma}

\begin{proof}
Suppose that $d\leqslant 2$. It follows from
Corollary~\ref{corollary:from-the-number-of-good-points} that
$$
2\geqslant d=2\mu-3-\sum_{i=k}^{l}ic_{i}\geqslant\mu-2\geqslant 1,
$$
but $\mu\ne 3$ by Lemma~\ref{lemma:swapping}, because
$|\Gamma|\leqslant|\Sigma\setminus\Lambda|\leqslant
\mu(\mu-\sum_{i=k}^{l}c_{i}i)-2$.

Thus, we see that $\mu=4$. Then $k=3$ by Lemma~\ref{lemma:conic},
which implies that
$$
\big|\Gamma\big|\leqslant\big|\Sigma\setminus\Lambda\big|\leqslant 14-4\sum_{i=k}^{l}c_{i}i\leqslant 2,%
$$
which is impossible by Lemma~\ref{lemma:swapping}, because
$d\geqslant 1$.
\end{proof}

It follows from the
inequality~\ref{equation:number-of-good-points} that
$|\Gamma^{\prime}|=|\Gamma|\leqslant|\Sigma\setminus\Lambda|\leqslant\mu(\mu-\sum_{i=k}^{l}c_{i}i)-2$.
Then
$$
\big|\Gamma^{\prime}\big|\leqslant\mu\left(\mu-\sum_{i=k}^{l}c_{i}i\right)-2\leqslant\mathrm{max}\left\{\Big\lfloor\frac{d+3}{2}\Big\rfloor\left(d+3-\Big\lfloor\frac{d+3}{2}\Big\rfloor\right)-1,\Big\lfloor\frac{d+3}{2}\Big\rfloor^{2}\right\},%
$$
because $d=2\mu-3-\sum_{i=k}^{l}c_{i}i$ and $\mu\geqslant 3$.

\begin{lemma}
\label{lemma:lines} At most $d$ points of the set $\Gamma$ is
contained in a~line.
\end{lemma}

\begin{proof}
Suppose that at least $d+1$ points of the set $\Gamma$ is
contained in some~line. Then
$$
\mu\geqslant d+1=2\mu-2-\sum_{i=k}^{l}c_{i}i,
$$
because at most $\mu$ points of $\Gamma$ is contained in a~line.
It follows from
Corollary~\ref{corollary:from-the-number-of-good-points} that
$$
\mu-1\geqslant\sum_{i=k}^{l}c_{i}i\geqslant \mu-2.
$$

Suppose that $\sum_{i=k}^{l}c_{i}i=\mu-2$. Then $|\Gamma|\leqslant
2\mu-2$. So, the set $\Gamma$ imposes independent linear
conditions on forms of degree $d=\mu-1$ by Theorem~2 in
\cite{EiJ87}, which is impossible by~Lemma~\ref{lemma:swapping}.

We see that $\sum_{i=k}^{l}c_{i}i=\mu-1$. Then $|\Gamma|\leqslant
\mu-2=d$, which is impossible by Lemma~\ref{lemma:swapping}.
\end{proof}

Therefore, at most $d$ points of the set $\Gamma^{\prime}$ lies on
a line by
Lemmas~\ref{lemma:lines}~and~\ref{lemma:zero-dimensional}.

\begin{lemma}
\label{lemma:final} For every $t\leqslant (d+3)/2$, at most
$$
t\big(d+3-t\big)-2
$$
points of the set $\Gamma^{\prime}$ lie on a curve  of degree $t$
in $\Pi\cong\mathbb{P}^{2}$.
\end{lemma}

\begin{proof}
At most $\mu t$ points of the set $\Gamma^{\prime}$ lie on a curve
of degree $t$. It is enough to show that
$$
t\big(d+3-t\big)-2\geqslant \mu t
$$
for every $t\leqslant (d+3)/2$ such that $t>1$ and
$t(d+3-t)-2<|\Gamma^{\prime}|$. But
$$
t(d+3-t)-2\geqslant t\mu\iff \mu-\sum_{i=k}^{l}c_{i}i>t,
$$
because $t>1$. Thus, we may assume that
$t(d+3-t)-2<|\Gamma^{\prime}|$ and
$$
\mu-\sum_{i=k}^{l}c_{i}i\leqslant t\leqslant {\frac{d+3}{2}}.%
$$

Let $g(x)=x(d+3-x)-2$. Then
$$
g\big(t\big)\geqslant g\Big(\mu-\sum_{i=k}^{l}c_{i}i\Big),
$$
because $g(x)$ is increasing for $x<(d+3)/2$. Therefore, we have
$$
\mu\Big(\mu-\sum_{i=k}^{l}ic_{i}\Big)-2\geqslant\big|\Gamma^{\prime}\big|>g(t)\geqslant
g\Big(\mu-\sum_{i=k}^{l}c_{i}i\Big)=\mu\Big(\mu-\sum_{i=k}^{l}ic_{i}\Big)-2,
$$
which is a contradiction.
\end{proof}

Thus, the set $\Gamma^{\prime}$ imposes independent linear
conditions on forms of degree $d$ by Theorem~\ref{theorem:Bese},
which implies that the set $\Gamma$ also imposes independent
linear conditions on forms of degree~$d$, which is impossible by
Lemma~\ref{lemma:swapping}. The assertion of
Theorem~\ref{theorem:main} is proved.

\end{document}